\newenvironment{E}{\begin{equation}}{\end{equation}}
\def\proof{\noindent{\bf Proof: }}
\def\qed{ \hskip 20pt{\vrule height7pt width6pt depth0pt}\hfil}
\def\forb{{\hbox{forb}}}
\def\fb{{\hbox{forb}}}
\def\forbmax{{\hbox{forbmax}}}
\def\fbmax{{\hbox{forbmax}}}

\def\0{{\bf 0}}
\def\1{{\bf 1}}

\def\Av{{\mathrm{Avoid}}}

\def\BB{{\mathrm{BB}}}

\newcommand{\linelessfrac}[2]{\genfrac{}{}{0pt}{}{#1}{#2}}
\newcommand{\ncols}[1]{\| #1 \|}
\newcommand{\rf}[1]{(\ref{#1})}
\newcommand{\trf}[1]{Theorem~\ref{#1}}
\newcommand{\lrf}[1]{Lemma~\ref{#1}}
\newcommand{\crf}[1]{Corollary~\ref{#1}}

\newcommand{\srf}[1]{Section~\ref{#1}}

\newcommand{\probrf}[1]{Problem~\ref{#1}}
\documentclass[12pt]{article}
\newtheorem{thm}{Theorem}[section]
\newtheorem{lemma}[thm]{Lemma}

\newtheorem{cor}[thm]{Corollary}

\newtheorem{probl}[thm]{Problem}

\usepackage{amssymb}
\usepackage{amsmath}
\title{ Multivalued Matrices and Forbidden Configurations }
\author{R.P. Anstee\thanks{Research supported in part by
NSERC} , Jeffrey Dawson \thanks{Research supported in part by
NSERC USRA}
\\Mathematics Department\\The University of British Columbia\\Vancouver,
B.C. Canada V6T 1Z2\\ \\
Linyuan Lu\\University of South Carolina\\Columbia SC, USA  \\ \\
Attila Sali\\R\'enyi Institute\\Budapest, Hungary  \\ }

\begin{document}
\maketitle
\begin{abstract}
An $r$-matrix is a matrix with symbols in $\{0,1,\ldots,r-1\}$. A matrix is simple if it has no repeated columns. Let ${\cal F}$ be a finite set of $r$-matrices. Let $\forb(m,r,{\cal F})$ denote the maximum number of columns possible in a simple $r$-matrix
 $A$ that has no submatrix which is a row and column permutation of any $F\in{\cal F}$.  Many investigations have involved $r=2$. For  general $r$, $\forb(m,r,{\cal F})$ is polynomial in $m$ if and only if for every pair $i,j\in\{0,1,\ldots,r-1\}$ there is a matrix in ${\cal F}$ whose entries are only $i$ or $j$.  
 Let ${\cal T}_{\ell}(r)$ denote the following $r$-matrices. For a pair $i,j\in\{0,1,\ldots,r-1\}$ we form four $\ell\times\ell$ matrices namely the matrix with $i$'s on the diagonal and $j$'s off the diagonal and the matrix with $i$'s on and above the diagonal and $j$'s below the diagonal and the two matrices with the roles of $i,j$ reversed.
Anstee and Lu determined that $\forb(m,r,{\cal T}_{\ell}(r))$ is a constant.  Let ${\cal F}$ be a finite set of 2-matrices. We ask if
$\forb(m,r,{\cal T}_{\ell}(3)\backslash {\cal T}_{\ell}(2)\cup {\cal F})$ is 
$\Theta(\forb(m,2,{\cal F}))$ and settle this in the affirmative for some cases including most 2-columned $F$.

\end{abstract}

\section{Introduction}

We define 
 a matrix to be \emph{simple} if it has no repeated columns. A (0,1)-matrix that is simple is the matrix analogue of a set system (or simple hypergraph)
 thinking of the matrix as the element-set incidence matrix.  
 We generalize to allow more entries in our matrices and define
 an $r$-\emph{matrix} be a matrix whose entries are in $\{0,1,\ldots, r-1\}$. We can think of this as an $r$-coloured matrix. 
 For $r=2$, $r$-matrices are (0,1)-matrices and for $r=3$, $r$-matrices are (0,1,2)-matrices.  We examine extremal problems and let $\ncols{A}$ denote the number of columns in $A$.

We will use the language of matrices in this paper rather than sets. For two  matrices $F$ and $A$, we write $F \prec A$, and say that $A$ has $F$ as a \emph{configuration}, if there is a  submatrix of $A$ which is a row and column permutation of $F$. Row and column order matter to submatrices but not to configurations.  Let ${\cal F}$ denote a finite set of matrices. Let 
$$\Av(m,r,{\cal F})=\left\{A\,:\,A\hbox{ is }m\hbox{-rowed and simple }r\hbox{-matrix}, F\not\prec A\hbox{ for }F\in{\cal F}\right\}.$$
Our extremal function of interest is
$$\forb(m,r,{\cal F})=\max_A\{\ncols{A}\,:\,A\in\Av(m,r,{\cal F})\}.$$
In the case $r=2$, we are considering (0,1)-matrices and then we drop $r$ from the notation to write
$\Av(m,2,{\cal F})=\Av(m,{\cal F})$ and $\forb(m,2,{\cal F})=\forb(m,{\cal F})$.
We define 
$$\forbmax(m,r,{\cal F})=\max_{m'\le m}\forb(m',r,{\cal F}).$$
 It has been conjectured by Anstee and Raggi 
\cite{miguelthesis} that $\forbmax(m,2,{\cal F})=\forb(m,2,{\cal F})$ for large $m$ (which is a type of monotonicity). For many ${\cal F}$ this is readily proven.

The following dichotomy between polynomial and exponential bounds is striking. Denote an $(i,j)$-matrix as a matrix whose
 entries are $i$ or $j$.

\begin{thm} (F\"uredi and Sali \cite{FS})  Let ${\cal F}$ be a family of $r$-matrices.  If  for every pair $i,j\in\{0,1,\ldots ,r-1\}$, there is an $(i,j)$-matrix in ${\cal F}$ then for some $k$, $\forb(m,r,{\cal F})$ is 
$O(m^k)$. If there is some pair $i,j\in\{0,1,\ldots ,r-1\}$ so that ${\cal F}$ has no $(i,j)$-matrix then $\forb(m,r,{\cal F})$ is 
$\Omega(2^m)$.\label{dichotomy}\end{thm}

It would be of interest to have more examples of forbidden families of configurations where we can determine the asymptotics of
$\forb(m,r,{\cal F})$. There are known examples given in \cite{FS}.   There is a generalization of a result of Balogh and Bollob\'as \cite{BaBo} for (0,1)-matrices to $r$-matrices.  Define the generalized identity matrix
 $I_{\ell}(a,b)$ as the $\ell\times\ell$ $\,\,
r$-matrix with $a$'s on the diagonal and $b$'s elsewhere. The standard identity matrix is $I_{\ell}(1,0)$. 
Define the generalized triangular matrix 
$T_{\ell}(a,b)$ as the $\ell\times\ell$ 
$\,\,r$-matrix with $a$'s below the diagonal and $b$'s elsewhere.  The standard upper triangular matrix is $T_{\ell}(0,1)$.
Let
$${\cal T}_{\ell}(r)=\left\{ I_{\ell}(a,b)\,:\,a,b\in
  \{0,1,\cdots, r-1\}, a\ne b \right\}$$
  \begin{E}\bigcup 
  \left\{ T_{\ell}(a,b)\,:\,a,b\in
  \{0,1,\cdots, r-1\}, a\ne b \right\}.
\label{defn}\end{E}
By \trf{dichotomy},  $\forb(m,r,{\cal T}_{\ell}(r))$ is bounded by a polynomial but much more is true.

\begin{thm}\cite{ALu} Given $r,\ell$, there is a constant $c(r,\ell)$ so that $\forb(m,r,{\cal T}_{\ell}(r))\le c(r,\ell)$.
\label{constantbd}\end{thm}
We will use the constant $c(r,\ell)$ repeatedly in this paper. This is a kind of Ramsey Theorem, a particular structured configuration appears in any $r$-matrix of a suitably large number of distinct columns.  An important result is that $c(r,\ell)$ is $O(2^{c_r\ell^2})$ for some constant $c_r$.  Not unexpectedly, Ramsey Theory shows up in the proof.  \srf{ramsey} contains a number of proofs using  Ramsey theory.

${\cal T}_{\ell}(2)$ consists of (0,1)-matrices (i.e. 2-matrices).  This paper considers forbidding the matrices ${\cal T}_{\ell}(r)\backslash {\cal T}_{\ell}(2)$. Note that  any (0,1)-matrix $A\in\Av(m,r,{\cal T}_{\ell}(r)\backslash {\cal T}_{\ell}(2))$ and so
$\forb(m,r,{\cal T}_{\ell}(r)\backslash {\cal T}_{\ell}(2))=\Omega(2^m)$.  Forbidding ${\cal T}_{\ell}(r)\backslash {\cal T}_{\ell}(2)$ may be somewhat like  asking the matrices to be (0,1)-matrices.  
\begin{thm} Let $r,\ell$ be given. Then $\forb(m,r,{\cal T}_{\ell}(r)\backslash {\cal T}_{\ell}(2))$ is $\Theta(2^m)$.\label{noF}\end{thm}
\proof A construction in $\Av(m,r,{\cal T}_{\ell}(r)\backslash {\cal T}_{\ell}(2))$  is to take all $2$-columns on $m$ rows. 

Take any matrix $A\in\Av(m,r,{\cal T}_{\ell}(r)\backslash {\cal T}_{\ell}(2))$ and replace all entries $2,3,\ldots ,r-1$ by $1$'s  to obtain the $2$-matrix $A'$, not necessarily simple.  The number of different columns in $A'$ is at most $2^m$. 

Let $\alpha$ be a column of $A'$. 
Let $B$ denote the submatrix of $A$ consisting of all columns of $A$ that map to $\alpha$ under the replacements.  Let $B'$ be the simple submatrix of $B$ consisting of the rows of $B$ where $\alpha$ has 1's.  Then $\ncols{B}=\ncols{B'}\le c(r-1,\ell)$ else we have a configuration in ${\cal T}_{\ell}(r)\backslash {\cal T}_{\ell}(2)$ in $B'$ (using \trf{constantbd}) with symbols chosen from 
$\{1,2,\ldots,r-1\}$. 

Combining these two observations yields the desired bound.\qed

\vskip 10pt
By the same argument we can show $\forb(m,r,{\cal T}_{\ell}(r)\backslash {\cal T}_{\ell}(s))$ is $\Theta(s^m)$ but 
the focus is on $s=2$ in this paper.  In this paper we will also take $r=3$. Note that \lrf{rto3} provides a justification for this restriction. Define the matrices $T_{\ell}(a,b,c)$ as the $\ell\times\ell$ matrix with $a$'s below the diagonal, $b$'s on the diagonal and $c$'s above the diagonal. In our problems we can require $a\ne b$. These appear in the proof of \trf{constantbd} but, for $a\ne b\ne c$,  are not matrices of just two entries which are referred to in \trf{dichotomy}.  One general result in this direction is the following.  

\begin{thm}\cite{ALu15} Let ${\cal F}$ be a finite family of (0,1)-matrices. Then  \hfill\break$\forb(m,3,{\cal T}_{\ell}(3)\backslash {\cal T}_{\ell}(2)\cup T_{\ell}(0,2,1)\cup {\cal F})$ 
is $O(\forbmax(m,{\cal F}))$.
\label{extra}\end{thm}

Another version of \trf{extra} with restricted column sums (\emph{column sum} will refer in this setting to the number of 1's)  is given in \srf{ramsey} with the analogous proof. We are not pleased with the inclusion of $T_{\ell}(0,2,1)$ in \trf{extra} and think it can be avoided. 

\begin{probl} Let $F$ be a (0,1)-matrix. Is it true that $\forb(m,3,{\cal T}_{\ell}(3)\backslash {\cal T}_{\ell}(2)\cup F)$ is $\Theta(\forbmax(m,F))$?
\label{problem}\end{probl}

Obviously the configuration $T_{\ell}(0,2,1)$ will be problematical. We will let $\ell$ take on large but constant values. Some results given below  support a yes answer. For example if we forbid nothing in the (0,1)-world then the maximum number of possible distinct (0,1)-columns is $2^m$.  One could say that
``$\forbmax(m,\emptyset)=2^m\,\,$''. Now using \trf{noF}, we see that \probrf{problem} is true in this case.

  Given $s=2$, one can show it suffices to consider  $r=3$ in \probrf{problem}. The  argument is similar to \trf{noF} and uses Ramsey Theory. The proof is given in \srf{ramsey}. 
\begin{lemma}  Let $r>2$ and  $\ell$ be given. Then there is a constant $bd(\ell)$ so that \hfil\break$\forb(m,r,{\cal T}_{\ell}(r)\backslash {\cal T}_{\ell}(2)\cup F)$ is
$O(\forb(m,3,{\cal T}_{bd(\ell)}(3)\backslash {\cal T}_{bd(\ell)}(2)\cup F))$.  \label{rto3}\end{lemma}

Given the answer `yes' to \probrf{problem}, this yields a justification for restricting to $r=3$.  The argument could also be extended to ${\cal T}_{\ell}(r)\backslash {\cal T}_{\ell}(s)$ but 
the focus is on $s=2$.
 Many configurations $F$ can be handled by \trf{intriangle} and in particular configurations with more than two columns.  

\begin{thm}Let $F\prec T_{\ell/2}(0,1)$. Then $\forb(m,3,{\cal T}_{\ell}(3)\backslash {\cal T}_{\ell}(2)\cup F)$ is $\Theta(\forbmax(m,F))$.
\label{intriangle}\end{thm}
\proof Note that $T_{\ell/2}(0,1)\prec T_{\ell}(0,2,1)$ by considering the submatrix of $T_{\ell}(0,2,1)$ consisting of the  even indexed columns and the odd indexed  rows. Thus if $F\not\prec A$, then $T_{\ell}(0,2,1)\not\prec A$. Apply \trf{extra}. \qed
\vskip 10pt

One important corollary is the following. 
\begin{cor}$\forb(m,3,{\cal T}_{\ell}(3)\backslash {\cal T}_{\ell}(2)\cup [0\,1])$ is $\Theta(1)$. \label{[01]}\end{cor}
\proof $[0\,1]\prec T_{\ell/2}(0,1)$ for $\ell\ge 4$. \qed
\vskip 10pt

This paper provides a number more results in this direction mostly involving configurations of two columns.
Define $F_{a,b,c,d}$ to be the $(a+b+c+d)\times 2$ configuration with $a$ rows $[1\,1]$, $b$  rows $[1\,0]$, 
$c$  rows $[0\,1]$, and $d$  rows $[0\,0]$. The asymptotics of $\forb(m,F_{a,b,c,d})$ have been completely determined by Anstee and Keevash  \cite{AK}. Note that we can assume $a\ge d$ since otherwise we can tale the (0,1)-complement $F_{a,b,c,d}^c=F_{d,c,b,a}$. Also we may assume $b\ge c$ since as configurations $F_{a,b,c,d}=F_{a,c,b,d}$.  We note that 
$\forb(m,F_{a,b,0,0})$ is $\Omega(m^{a+b-1})$  by taking all columns of column sum $a+b$  and  a different construction shows 
$\forb(m,F_{0,b,b,0})$ is 
$\Omega(m^b)$.     The important upper bounds  are for $a\ge 1$, $\forb(m,F_{a,b,b,a})$ is $\Theta(m^{a+b-1})$ \cite{AK} and 
$\forb(m,F_{0,b+1,b,0})$ is $\Theta(m^b)$ \cite{AK}. 
Note that $I_2=F_{0,1,1,0}$. This is the first result not covered by \trf{intriangle}.

\begin{thm}$\forb(m,3,{\cal T}_{\ell}(3)\backslash {\cal T}_{\ell}(2)\cup I_2)$ is $\Theta(\forbmax(m,I_2))$. 
\label{I2}\end{thm}

\begin{thm}Let $a\ge 0$ and $b\ge 2$ be given. Then $\forb(m,3,{\cal T}_{\ell}(3)\backslash {\cal T}_{\ell}(2)\cup F_{a,b,b,a})$ is 
$\Theta(\forbmax(m,F_{a,b,b,a}))$. \label{abba}\end{thm}

We give the proofs in \srf{twocolumn}.   Note the subtlety that  $\forbmax(m,F_{0,b,b,0})$
is 
$\Theta(m^b)$ where as, for $a\ge 1$,  $\forbmax(m,F_{a,b,b,a})$ is 
$\Theta(m^{a+b-1})$. The proofs use  results for two columned forbidden configurations from \cite{AK}.
The other critical two columned result concerns $F=F_{0,b+1,b,0}$ for which we don't know the answer for \probrf{problem}.

Define $t\cdot F=[F\,F\,\cdots F]$ to be the concatenation of $t$ copies of $F$.

\begin{thm} Let $F$ be a given $k\times p$ (0,1)-matrix. 
Then $\forb(m,3,{\cal T}_{\ell}(3)\backslash {\cal T}_{\ell}(2)\cup t\cdot F)$ is $O(\max\{m^{k},\forb(m,3,{\cal T}_{\ell}(3)\backslash {\cal T}_{\ell}(2)\cup F)\})$. \label{ttimes}\end{thm}
\proof  Let $A\in\Av(m,3,{\cal T}_{\ell}(3)\backslash {\cal T}_{\ell}(2))$ with 
$$\ncols{A}>
(t-1)p\binom{m}{k}+\forb(m,3,{\cal T}_{\ell}(3)\backslash {\cal T}_{\ell}(2)\cup F)+1.$$
 Then 
$F\prec A$. Remove from $A$ the $p$ columns containing a copy of $F$ and repeat.  We will generate at least
$(t-1)\binom{m}{k}+1$ copies of $F$ and hence at least $t$ column disjoint copies of $F$ in the same set of $k$ rows and so $t\cdot F\prec A$. \qed

\vskip 10pt
To apply this, we need to know  $\forb(m,3,{\cal T}_{\ell}(3)\backslash {\cal T}_{\ell}(2)\cup F)$. The following is established in 
\srf{twocolumn}.

\begin{thm} Let 
\begin{E}H=\left[\begin{array}{ccc}0&1&0\\ 0&0&1\\ \end{array}\right]\label{H}\end{E}
Then $\forb(m,3,{\cal T}_{\ell}(3)\backslash {\cal T}_{\ell}(2)\cup H)$ is $\Theta(m)$.
\label{linear}\end{thm}

\begin{cor} Given $H$ in \rf{H}, we have 
$\forb(m,3,{\cal T}_{\ell}(3)\backslash {\cal T}_{\ell}(2)\cup t\cdot H)$ is $\Theta(m^2)$.\label{ttimeslinear}\end{cor}
\proof We apply \trf{ttimes} and \trf{linear}. \qed
\vskip 10pt
\trf{linear} and \crf{ttimeslinear} are yes instances of \probrf{problem} since $\forb(m,H)$ is $\Theta(m)$ and
$\forb(m,t\cdot H)$ is $\Theta(m^2)$ \cite{survey}.  

Given an $m_1\times n_1$ matrix $A$ and a $m_2\times n_2$ matrix $B$,  defne the product of two matrices $A\times B$ as the 
$(m_1+m_2)\times n_1n_2$ matrix obtained from placing each column of $A$ on top of each column of $B$ 
for all possible pairs of columns. Let $F$ be given with
$$
0\times 1\times F=\left[\begin{array}{c}0\,0\,\cdots 0\\ 1\,1\,\cdots 1\\  F\\ \end{array}\right]$$
In \cite{survey}, we establish that $\forb(m,0\times 1\times F)$ is $O(m\cdot\forb(m,F))$. We establish this version of the \probrf{problem} in \srf{0x1xF}.

\begin{thm}
$\forb(m,3,{\cal T}_{\ell}(3)\backslash {\cal T}_{\ell}(2)\cup  0\times 1\times F)$ is
$O(m\cdot \forb(m,{\cal T}_{\ell}(3)\backslash {\cal T}_{\ell}(2)\cup  F))$ \label{rows0and1}\end{thm}
This result extends results for $F_{a,b,b,a}$ to $F_{a+1,b,b,a+1}$ and can be used in other instances such as $H$ above. We finish the paper with some open problems.

\section{Results using Ramsey Theory}\label{ramsey}

We apply Ramsey Theory to help us find configurations in ${\cal T}_{\ell}(3)\backslash {\cal T}_{\ell}(2)$ etc.
We  use the $p$ colour Ramsey number $R_p({t_1,t_2,\ldots ,t_p})$ as the smallest number $n$ such that for every edge colouring of $K_n$ with $p$ colours there is some colour $i$ so that there is a clique of size $t_i$ with all edges of  colour $i$. Typical notation is that for $t_1=t_2=\cdots =t_p=t$, we write $R_p({t_1,t_2,\ldots ,t_r})=R_p(t^p)$.  While these  numbers can be  large, we can for example bound $R_p(t^p)\le 2^{pt}$.

Let $r,s$ be given integers with $r>s\ge 2$. Let us define a set ${\cal P}_{t}^x(r)$ of $t\times t$ matrices by the following template which will have choices $x,y_1,y_2,\ldots ,y_{t}\in\{1,2,\ldots,r-1\}$ where we require $y_j\ne x$ for $j\in[t]$.  The entries marked * may be given entries in $\{0,1,\ldots ,r-1\}$ in any possible way.

\begin{E}
{\cal P}^x_{t}:\quad\left[\begin{array}{cccccc}
y_1&&&&&\\
x&y_2&&&*&\\
x&x&y_3&&&\\
\vdots&&&\ddots&\\
x&x&x&&y_{t -1}&\\
x&x&x&\cdots&x&y_{t}\end{array}\right]\label{Ptemplate}\end{E}

\begin{lemma} Let $\ell,r,s$ be given with $r>s\ge 2$.  Let $t=(r-1)(R_r((2\ell)^r)-1)-1$. Assume $A$ is an $m$-rowed simple $r$-matrix. Assume there is some $G\in {\cal P}_{t}^x$ with $G\prec A$ and such that if $x\in\{0,1,\ldots ,s-1\}$ then $y_j\in\{s,s+1,\ldots ,r-1\}$ for all $j\in [t]$. 
Then there is some  $F$ with  $F\prec G$ and
$$F\in\left({\cal T}_{\ell}(r)\backslash{\cal T}_{\ell}(s)\right)\bigcup \{T_{\ell}(x,z,u)\, ;\,x,u\in \{0,1,\ldots ,s-1\}, x\ne u, z\notin 
\{0,1,\ldots ,s-1\}\}. $$ .
\label{newtriangular}\end{lemma}
\proof Assume there is some $G\in {\cal P}_{t}^x$ with $G\prec A$ and such that if $x\in\{0,1,\ldots ,s-1\}$ then $y_i\in\{s,s+1,\ldots ,r-1\}$ for all $i\in [t]$. 

 First  assume $x\notin\{0,1,\ldots ,s-1\}$. There are $r-1$ choices for each $y_j$ and hence there is some choice $z\in \{0,1,\ldots ,s-1\}\backslash x$ which appears at least
 $R_r((2\ell)^r)$ times on the diagonal. Now form a graph whose vertices are the rows $i$ with $y_j=z$ and we colour edge $a,b$ for $a<b$ by the entry in the $a,b$ location of $G$ (above the diagonal). There will be at least $R_r((2\ell)^r)$ vertices and there will be at most $r$ colours and so by the Ramsey number  there will be a clique of size $2\ell$ of all edges of the same colour, say colour $u$. If $u=z$ we have
  $T_{2\ell}(x,z)\prec A$. If $u=x$ we have $I_{2\ell}(x,z)\prec A$. If $u\ne x,z$ then we consider the configuration $T_{2\ell}(x,u,z)$ of size $2\ell$ induced by the clique and  the even columns and the odd rows to show $T_{\ell}(x,y)\prec A$. All three cases yield a configuration 
  in ${\cal T}_{\ell}(r)\backslash{\cal T}_{\ell}(s)$.
  
  Now assume $x\in\{0,1,\ldots ,s-1\}$ then there are $r-s-1$ choices for each $y_j$ and hence there is some choice $z\notin \{0,1,\ldots ,s-1\}$ which appears at least
 $R_r((2\ell)^r)$ times on the diagonal. Now we proceed as above to obtain a configuration $T_{2\ell}(x,z,u)$.  If $u\in\{s,s+1,\ldots ,r-1\}$ then we obtain a configuration in ${\cal T}_{\ell}(r)\backslash{\cal T}_{\ell}(s)$. If $u=x$ we obtain a configuration $I_{2\ell}(x,z)$ which is 
 in ${\cal T}_{\ell}(r)\backslash{\cal T}_{\ell}(s)$.  If 
 $u\in\{0,1,\ldots ,s-1\}$ with $u\ne x$, then we obtain a configuraton $T_{2\ell}(x,z,u)$ with $x,u\in \{0,1,\ldots ,s-1\}$ and $x\ne u$ (which does not yield a configuration in
 ${\cal T}_{\ell}(r)\backslash{\cal T}_{\ell}(s)$). 
 \qed

\vskip 10pt

Our application of the \lrf{newtriangular} to \trf{extra} will be in the case $r=3$ and $s=2$ and then
$ \{T_{\ell}(x,z,u)\, ;\,x,u\in \{0,1,\ldots ,s-1\}, x\ne u, z\notin 
\{0,1,\ldots ,s-1\}$ is the single configuration $T_{\ell}(0,2,1)$.  We prove in greater generality.

\vskip 10pt

\noindent{\bf Proof of \trf{extra}}:\hskip .5cm 
The idea of the proof is to use the induction to generate configurations corresponding to 
matrices in ${\cal P}^x_{t}$  that enable us to apply the proof of \lrf{newtriangular} and obtain 
 matrices in ${\cal T}_{\ell}(r)\backslash {\cal T}_{\ell}(s)$.

We use the following function $f$ in our proof.  Let $f$ be determined by the recurrence 
\begin{E}f(p_0,p_1,\ldots ,p_{r-1})=\sum_{i=0}^{r-1}f(p_0,p_1,\ldots,p_i-1,\ldots, p_{r-1}),\label{recurrence}\end{E}
and the base cases that $f(p_0,p_1,\ldots ,p_{r-1})=1$ if $p_i=1$ for any $i\in\{0,1,\ldots ,r-1\}$.  Solving this  exactly seems difficult but since $f$ satisfies the same recurrence as the multinomial coefficients, with smaller base cases, we obtain
\begin{E}f(p_0,p_1,\ldots ,p_{r-1})\le \frac{(p_0+p_1+\cdot +p_{r-1}-r)!}{(p_0-1)!(p_1-1)!\cdots(p_{r-1}-1)!}\label{formula}\end{E}
Let $g(p_0,p_1,\ldots ,p_{r-1})=$  $f(p_0,p_1,\ldots ,p_{r-1})\cdot\forbmax(m,{\cal F})$. 

We will establish for fixed $m$ but by induction on $\sum_ip_i$, that if $A$ is an $n$-rowed simple $r$-matrix with $n\le m$ and $\ncols{A}> g(p_0,p_1,\ldots ,p_{r-1})$ then for some $i\in\{0,1,\ldots ,r-1\}$, $A$ will contain configuration $F\in{\cal F}$ or a configuration in 
${\cal P}^i_{p_i}$ satisfying the condition that if $i\in\{0,1,\ldots ,s-1\}$, then $y_j\in\{s,s+1,\ldots ,r-1\}$ for $j\in [P_i]$. 
We use $\forbmax$ so that $\forbmax(m,s,{\cal F})\ge \forb(n,s,{\cal F})$.

If $p_i=1$, then an element of  ${\cal P}^i_{p_i}$ is a $1\times 1$ matrix. For $i\in\{0,1,\ldots ,s-1\}$, then the entry in the $1\times 1$ matrix must not be in 
$\{0,1,\ldots ,s-1\}$ and if $i\notin\{0,1,\ldots ,s-1\}$, then the entry in the $1\times 1$ matrix must not be $i$. In the former case, we require the matrix to have some entry not in
$\{0,1,\ldots ,s-1\}$ which would only be difficult if $A$ was an $s$-matrix. In that case $\ncols{A}\le\forb(n,s,{\cal F})\le\forbmax(m,s,{\cal F})$ and we note that $f(p_0,p_1,\ldots ,p_{r-1})=1$ for $p_i=1$.  In the latter case we are merely requiring that the matrix $A$ has at least two different entries which would only not occur for $\ncols{A}=1$. In either case we are able to obtain an instance of ${\cal P}_1^i$  in $A$ if $\ncols{A}> g(p_0,p_1,\ldots ,p_{r-1})$. This establishes the required base cases for the induction.

  Assume $p_i\ge 2$ or all $i\in\{0,1,\ldots, r-1\}$. Consider a matrix 
$A\in\Av(n,r,{\cal P}_{p_0}^0\cup{\cal P}_{p_1}^1\cup\cdots\cup{\cal P}_{p_{r-1}}^{r-1}\cup {\cal F})$ with $n\le m$ and 
$\ncols{A}>g(p_0,p_1,\ldots ,p_{r-1})$. We wish to obtain a contradiction. 

Choose a row $w$ of $A$ which has at least two different entries one of which is not in 
$\{0,1,\ldots ,s-1\}$.  If there is no such row then either $\ncols{A}=1$ or $A$ is an $s$-matrix.  In the latter case, we have $\ncols{A}
> g(p_1,p_2,\ldots p_{r-1})\ge \forbmax(m,s,{\cal F})\ge \forb(n,s,{\cal F})$ and so $F\prec  A$, a contradiction. We may assume a row $w$ of $A$, which has at least two different entries one of which is not in 
$\{0,1,\ldots ,s-1\}$, exists.

Decompose $A$ as follows by permuting rows and columns
\begin{E}A=\begin{array}{r}\hbox{row }w\rightarrow\\ \\\end{array}
\left[\begin{array}{c|c|c|c|c}
0\,0\,\cdots 0&1\,1\,\cdots 1&2\,2\,\cdots 2&\cdots&r-1\,r-1\,\cdots r-1\\
G_0&G_1&G_{2}&\cdots&G_{r-1}\\
\end{array}\right]
\label{wdecomp}\end{E}
 Each $G_i$ is simple. Now 
$$\ncols{A}
=\sum_{i=0}^{r-1}\ncols{G_i}
>g(p_0,p_1,\ldots,p_{r-1})=f(p_0,p_1,\ldots ,p_{r-1})\cdot\forbmax(m,s,{\cal F})$$
$$=\left(\sum_{i=0}^{r-1}f(p_0,p_1,\ldots,p_i-1,\ldots, p_{r-1})\right)\cdot\forbmax(m,s,{\cal F}).$$
From the recurrence \rf{recurrence}, there is  some  $i$  with 
$$\ncols{G_i}>  g(p_0,p_1,\ldots,p_i-1,\ldots,p_{r-2},p_{r-1}).$$
Certainly $G_i\prec A$ and $G\in\Av(n-1,3,{\cal P}_{p_0}^0\cup{\cal P}_{p_1}^1\cup\cdots\cup{\cal P}_{p_{r-1}}^{r-1}\cup {\cal F})$.  
 Then by induction on $\sum_ip_i$, we can assume $G_i$ has a copy of ${\cal P}^i_{p_i-1}$ using the template \rf{Ptemplate} with 
  $x=i$ and  if $x=i\in\{0,1,\ldots ,s-1\}$ then  
$y_j\in\{s,s+1,\ldots ,r-1\}$ for all $j=1,2,\ldots ,p_i-1$.    We can extend to a copy of   ${\cal P}_{p_i}^i$ in $A$ by adding  row $w$ to extend by a row of $i$'s and then extend by a column from some $G_j$ with $j\ne i$.  If $i\in\{0,1,\ldots ,s-1\}$, then we can extend to a copy of   ${\cal P}^i_{p_i}$ in $A$ by adding  row $w$ to extend by a row of $i$'s and then extend by a column from some $G_h$ with 
$h\in\{s,s+1,\ldots ,r-1\}$.  This is possible since we have assumed that row $w$ has at least two different entries one of which is not in 
$\{0,1,\ldots ,s-1\}$.  Now some matrix $G$  in the family ${\cal P}^i_{p_i}$ has $G\prec A$. 

Specializing to $p_0=p_1=\cdots =p_{r-1}=
(r-1)(R_r((2\ell)^r)-1$ and applying \lrf{newtriangular} yields that $G$ contains a configuration in
$\left({\cal T}_{\ell}(r)\backslash{\cal T}_{\ell}(s)\right)\bigcup \{T_{\ell}(x,z,u)\, ;\,x,u\in \{0,1,\ldots ,s-1\}, x\ne u, z\notin 
\{0,1,\ldots ,s-1\}\}$ and then specializing to $r=3$ and $s=2$ yields the result. 
 \qed

\vskip 10pt
It was convenient to consider general $r,s$ but we will focus on $r=3$ and $s=2$.  The proof of
\trf{extra}  can be adapted to considering fixed column sum i.e. columns with a fixed number of 1's. 
In the case of 3-matrices, we define the \emph{column sum} of a 3-column $\alpha$ to be the number of 1's present. When there are no 2's in 
$\alpha$, this is the usual column sum.   Define
$$\forb_k(m,3,{\cal F}	)=\max\{\ncols{A}\,:\,A\in\Av(m,3,{\cal F}), \hbox{ all columns in }A\hbox{ have }k \,\,1's\},$$
and define $\forbmax_k$ similarly. 	
There are ${\cal F}$ for which we can exploit information about $\forb_k(m,3,{\cal F})$, deducing some information from
$\forb_k(m,{\cal F})$.

\begin{thm} Let ${\cal F}$ be a finite set of $(0,1)$-matrices. Let $\ell$ be given.  Then there exists a constant $d$ so that
\begin{E}\forb_k(m,3,{\cal T}_{\ell}(3)\backslash {\cal T}_{\ell}(2)\cup T_{\ell}(0,2,1)\cup {\cal F})\hbox{ is } 
O\left(\sum_{j=k-d}^k\forbmax_j(m,{\cal F})\right)\label{sumcolumnsum}\end{E}
\label{columnsumkextra}\end{thm}

\proof   We will follow the proof of \trf{extra}  but note how columns sums are affected.  	
Let $g_k(p_0,p_1,p_2)=  
f(p_0,p_1,p_2)\cdot \forbmax_k(m,F))$. 

Consider a matrix 
$A\in\Av_k(n,3,{\cal P}_{p_0}^0\cup{\cal P}_{p_1}^1\cup{\cal P}_{p_{2}}^{2}\cup {\cal F})$ with $n\le m$ and  $k>p+2R_3((2\ell)^3)$ and 
$\ncols{A}>g_k(p_0,p_1,p_2)$. We wish to obtain a contradiction. 

  It is convenient to interpret the proof of \trf{extra} as growing a tree where each node is associated with a matrix with three associated parameters $(p_0,p_1,p_2)$ and has some fixed column sum $s$. We begin with a root node corresponding to a matrix $A$ with parameters $(p,p,p)$ where $p=2R_3(2\ell,2\ell,2\ell)$. Then the matrices $G_0,G_1,G_2 $ can be viewed as the children. Our recursive growth of the tree begins with a node corresponding matrix $B$ for which we decompose by some row $w$ with at least two entries one of which is 2.  If we can't decompose then either $\ncols{B}=1$ or $B$ is an (0,1)-matrix. 
  
  Assume each column of $B$ has $s$ 1's.  Decompose $B$ as follows by permuting rows and columns
\begin{E}B=\begin{array}{r}\hbox{row }w\rightarrow\\ \\ \end{array}
\left[\begin{array}{c|c|c}
0\,0\,\cdots 0&1\,1\,\cdots 1&2\,2\,\cdots 2\\
H_0&H_1&H_{2}\\
\end{array}\right]
\label{wdecomp3}\end{E}
 Each $H_i$ is simple. Given that each column in $B$ has $s$  1's then  for each column in $H_0$ and $H_2$ has $s$ 1's and each column in $H_1$ has $s-1$ 1's.   Thus the nodes of our tree correspond to matrices with fixed column sum. 
 
 We also need to keep track of the current triple  $(q_0,q_1,q_2)$ for each node.    Thus if $B$ has the triple $(q_0,q_1,q_2)$ then $G_0$ has triple $(q_0-1,q_1,q_2)$,  $G_1$ has triple $(q_0,q_1-1,q_2)$ and $G_2$ has triple $(q_0,q_1,q_2-1)$.
 We do not decompose $B$ if $q_0=1$ or $q_1=1$ of $q_2=1$. Otherwise the node corresponding to $B$ has children $G_0,G_1,G_2$ with the possibility that $\ncols{G_0}=0$ or $\ncols{G_1}=0$ in which case $B$ would only have two children.

 Given the decomposition \rf{wdecomp3}, then $\ncols{A}$ is the sum of $\ncols{B}$ over all leaves $B$ of the tree. The leaves of the tree which cannot be further decomposed  correspond to matrices $B$ with $\ncols{B}=1$ or $B$ is a (0,1)-matrix or $B$  where
  the three parameters $(q_0,q_1,q_2)$ have  either $q_0=1$ or $q_1=1$ of $q_2=1$.
 
 We deduce that the depth of the tree is at most $d=3p=6R_3(2\ell,2\ell,2\ell)$ with a branching factor of $3$ and so there are at most $3^d$ nodes in the tree which is a constant.   Also we have that each node corresponds to a matrix with constant column sum  $s\in\{k-d,k-d+1,\ldots, k\}$ which is a constant cardinality set. 
 
 Now continue growing the tree until no further growth is possible. If the process generates a node $B$ with $q_0=1$ or $q_1=1$ of $q_2=1$, then by the arguments of \trf{extra}, there will be some configuration in ${\cal T}_{\ell}(3)\backslash {\cal T}_{\ell}(2)\cup T_{\ell}(0,2,1)$ in $B$ and hence in $A$.      A leaf node is one which corresponds to some (0,1)-matrix $B$  with constant column sum 
 $s\in\{k-d,k-d+1,\ldots, k\}$ for which we deduce that $\ncols{B}\le\forbmax_s(m,F)$. 
   
 The bound \rf{sumcolumnsum} now follows with the inclusion of some large constants.  
 \qed

\vskip 10pt
We will apply this result to 2-columned $F$.
\vskip 10pt
 \noindent {\bf Proof of \lrf{rto3}:} We readily note that 
$\forb(m,r,{\cal T}_{\ell}(r)\backslash {\cal T}_{\ell}(2)\cup F)\ge $ \hfil \break
$\forb(m,3,{\cal T}_{\ell}(r)\backslash {\cal T}_{\ell}(2)\cup F)$ since
$\Av(m,3,{\cal T}_{\ell}(3)\backslash {\cal T}_{\ell}(2)\cup F)\subseteq $\hfil\break
$\Av(m,r,{\cal T}_{\ell}(r)\backslash {\cal T}_{\ell}(2)\cup F)$.

Let $bd(\ell)=R_{(r-2)(r-2)}((2\ell)^{(r-2)(r-2)})$ where we assume $bd(\ell)>(r-2)\ell$. Let  \hfil\break$A\in\Av(m,r,{\cal T}_{\ell}(r)\backslash {\cal T}_{\ell}(2)\cup F)$. Replace all entries 
$3,4,\ldots ,r-1$ by $2$'s  to obtain $A'$.  The number of different columns in $A'$ is at most 
$\forb(m,3,{\cal T}_{bd(\ell)}(3)\backslash {\cal T}_{bd(\ell)}(2)\cup F)$ for the following reason.  If $F\prec A'$, then $F\prec A$ so we may assume $F\not\prec A'$.  Let $A''$ be the matrix obtained from $A'$ by keeping exactly one copy of each column. If $\ncols{A''}>\forb(m,3,{\cal T}_{bd(\ell)}(3)\backslash {\cal T}_{bd(\ell)}(2)\cup F)$ then there is configuration $G\prec A''$ with $G\in{\cal T}_{bd(\ell)}(3)\backslash {\cal T}_{bd(\ell)}(2)$.  There are several cases. 

If $G$ is a generalized identity matrix say $I_{bd(\ell)}(1,2)$, then in $A$, we have a configuration which has entries in $\{2,3,\ldots,r-1\}$ on the diagonal and 1's off the diagonal.  Then there is some entry $q\in \{2,3,\ldots,r-1\}$ appearing $\lceil bd(\ell)/(r-2)\rceil\ge \ell$ times (using $bd(\ell)>(r-2)\ell$) and we obtain a principal submatrix
of $G$ (row and column indices given by the diagonal entries $q$) in ${\cal T}_{\lceil bd(\ell)/(r-2)\rceil}(r)\backslash {\cal T}_{\lceil bd(\ell)/(r-2)\rceil}(2)$ in $A$. 

If  $G$ is a generalized identity matrix say $I_{bd(\ell)}(2,1)$,  then in $A$, we have a configuration which has entries in $\{2,3,\ldots,r-1\}$ off the diagonal and 1's on the diagonal.  Now apply Ramsey Theory by colouring a graph on $bd(\ell)$ vertices with the colour of edge $(i,j)$ for $i<j$  being the 2-tuple $a_{i,j},a_{j,i}$. There are $(r-2)^2$ colours and so if $bd(\ell)>
R_{(r-2)(r-2)}((2\ell)^{(r-2)(r-2)})$, then there is a clique of colour $p,q$ of size $2\ell$ and so  $2\ell\times 2\ell$ configuration whose entries on the diagonal are 1's and above the diagonal are $p$ and whose entries below the diagonal are $q$.  If $p=q$, we have a configuration 
in ${\cal T}_{2\ell}(r)\backslash {\cal T}_{2\ell}(2)$. If $p\ne q$, then we form an $\ell\times\ell$ configuration with $p$'s above the diagonal and $q$'s below the diagonal  (by taking even indexed columns and odd indexed rows) which is a configuration in ${\cal T}_{\ell}(r)\backslash {\cal T}_{\ell}(2)$.   Similar arguments handle the remaining cases.  

To determine the maximum number  of columns of $A$ that map into a given $(0,1,2)$-column $\alpha$ in $A'$,  let $\alpha$ have $t$ 2's and then the columns mapping into $\alpha$ correspond to a $t$-rowed simple matrix with entries in $\{2,3,\ldots,r-1\}$.  If the number of columns is bigger that $c(r-2,\ell)$, then those columns contain a configuration in ${\cal T}_{\ell}(r)$ whose entries are in $\{2,3,\ldots,r-1\}$ and so the configuration is in ${\cal T}_{\ell}(r)\backslash {\cal T}_{\ell}(2)$.  We now deduce that
$\ncols{A}\le c(r-2,\ell)\times forb(m,3,{\cal T}_{bd(\ell)}(3)\backslash {\cal T}_{bd(\ell)}(2)\cup F)$ yielding our bound. 
\qed
 \section{$0\times 1\times F$}\label{0x1xF}

Let $A\in\Av(m,3,{\cal T}_{\ell}(3)\backslash {\cal T}_{\ell}(2)\cup  0\times 1\times F)$.   
If we can choose a  pair of rows $i,j$ so that there are 
$\forb(m-2,3,{\cal T}_{\ell}(3)\backslash {\cal T}_{\ell}(2)\cup  F)+1$ columns of $A$ 
which have  0's on row $i$ and 1's in row $j$, then we have $F\prec A$, a contradiction. 
\begin{lemma} Let $\epsilon>0$ be given. Let $A$ be an $m$-rowed simple 3-matrix with each column having both a 0 and a 1 and at least 
$\epsilon m$ entries  either 0 or 1. Assume
\begin{E}\ncols{A}>
2\cdot\forbmax(m-2,3,{\cal T}_{\ell}(3)\backslash {\cal T}_{\ell}(2)\cup  F)
\frac{\binom{m}{2}}{\epsilon m-1}.\label{01pairs}\end{E}
Then $0\times 1\times F\prec A$.\label{em0's1's}\end{lemma}
\proof 
  We note that a column of $m$ rows  that has $p$ 0's and $q$  1's will have $pq$ pairs of rows $i,j$ containing the 
configuration $\left[\linelessfrac{0}{1}\right]$.  For a given $p,q$, the minimum number of configurations  
$\left[\linelessfrac{0}{1}\right]$ is $p+q-1$ when for example there is one 1 and $p+q-1$ 0's.  
An $m$-rowed column with at 
least one 0 and at least one 1 and at least $\epsilon m$ entries that are 0 or 1 will have
at least $\epsilon m-1$ configurations  
$\left[\linelessfrac{0}{1}\right]$.  There are $2\binom{m}{2}$ choices for $i,j$ when considered as an ordered pair.

If \rf{01pairs} is valid then 
then there will be a pair of rows $i,j$ with more than \hfil\break$2\cdot\forbmax(m,3,{\cal T}_{\ell}(3)\backslash {\cal T}_{\ell}(2)\cup  F)$ columns with the configuration $\left[\linelessfrac{0}{1}\right]$.  Thus there will be a pair of rows $i,j$ with at least $\forb(m,{\cal T}_{\ell}(3)\backslash {\cal T}_{\ell}(2)\cup  F)+1$ columns all with the submatrix $\left[\linelessfrac{0}{1}\right]$ (or all the reverse).  Then we can form an $(m-2)\times (\forbmax(m,3,{\cal T}_{\ell}(3)\backslash {\cal T}_{\ell}(2)\cup  F)+1)$ simple matrix $A'$ that when extends by a row of 0's and a row of 1's is contained in $A$.   SInce 
$A'\in\Av(m, {\cal T}_{\ell}(3)\backslash {\cal T}_{\ell}(2))$, we deduce that $F\prec A'$ and then $0\times 1\times F\prec A$, 
as desired. \qed

\vskip 10pt
 
\noindent {\bf Proof  of \trf{rows0and1}:}
If we have many columns with few 0's and 1's then we will show we are able to find in $A$  a $c\times c$ configuration $G$ in ${\cal P}^2_{c}$ of $A$ as in \rf{kcols} and then can use \lrf{newtriangular}.

Let $A\in\Av(m,3,{\cal T}_{\ell}(3)\backslash {\cal T}_{\ell}(2)\cup  F)$.  There are at most $c(2,\ell)$ (0,2)-columns and at most
$c(2,\ell)$ (1,2)-columns.  Let $A'$ be the matrix obtained from $A$ by deleting (0,2)-columns and (1,2)-columns.

 Now each column in $A'$ has  at least one 0 and one 1.  Let 
 \begin{E}\epsilon=\frac{1}{4R(2\ell,2\ell,2\ell)}.\label{epsilon}\end{E}
 Delete from $A'$ any rows entirely of 2's to obtain a simple matrix \hfil\break$A''\in\Av(t,3,{\cal T}_{\ell}(3)\backslash {\cal T}_{\ell}(2)\cup  F)$
 where $t\le m$.
  Let $A_2$ denote those columns of $A''$ with at most $\epsilon t$ 0's and 1's and let $A_{01}$ denote those columns of $A''$ with more than  $\epsilon t$ 0's and 1's .

We select columns of $A_2$ in turn to form the pattern ${\cal P}^2_c$ in \rf{Ptemplate} with $c=2\cdot R((2\ell)^{3})$. We can begin with a column on $(1-\epsilon)t$ 2's. At the $k$th stage we have $k$ columns (selected in the  order  displayed) with
\begin{E}
\begin{array}{c@{}l}
\overbrace{\begin{array}{@{}ccccc@{}}
\ne 2&&&&\\
2&\ne 2&&&\\
2&2&\ne 2&&\\
2&2&\cdots&\ne2&\\
2&2&\cdots&2&\ne 2\\
2&2&\cdots&2&2\\
2&2&\cdots&2&2\\
2&2&\cdots&2&2\\
\end{array}}^k&
\begin{array}{@{}c}
\begin{array}{c}
\\
\\
\\
\\
\\
\end{array}\\
\left.
\begin{array}{@{}l}
\\
\\
\\
\end{array}\right\}\ge(1-k\epsilon)t\\
\end{array}
\end{array}
\label{kcols}\end{E}
where the final block of 2's in rows $S$ has $|S|\ge (1-k\epsilon)t$.  Any column of $A_2$ 
not already chosen has 2's in at least $(1-(k+1)\epsilon)t$ rows of $S$.  To proceed we need that $\ncols{A_2}\ge c=2R_3(2\ell,2\ell,2\ell)$ and  we require that $(1-k\epsilon)t\ge 1$ for 
$k+1\le c=2R(2\ell,2\ell,2\ell)$. Our choice of $\epsilon$ \rf{epsilon} ensures this.
$A_2$ has no rows of 2's and so a column with a 0 or 1 in rows 
$S$ can be used to extend \rf{kcols} to the situation with $k+1$ columns.  We repeat until  we have $c=2R_3(2\ell,2\ell,2\ell)$ columns. Applying \lrf{newtriangular}, we obtain a matrix $F\in {\cal T}_{\ell}(3)$ with $F\prec A$ that has 2's below the diagonal and so we have obtained a configuration in ${\cal T}_{\ell}(3)\backslash{\cal T}_{\ell}(2)$, a contradiction.  
Thus 
\begin{E}\ncols{A_2}\le  2\cdot R_3(2\ell,2\ell,2\ell). \label{A_2bd}\end{E}

If 
$$\ncols{A_{01}}>
2\cdot\forbmax(m-2,3,{\cal T}_{\ell}(3)\backslash {\cal T}_{\ell}(2)\cup  F)
\frac{\binom{t}{2}}{\epsilon t-1}$$
then by \lrf{em0's1's}, $0\times 1\times F\prec A_{01}$. 
 Thus 
 \begin{E}\ncols{A_{01}}\le 
2\cdot\forb(m-2,{\cal T}_{\ell}(3)\backslash {\cal T}_{\ell}(2)\cup  F)
\frac{\binom{t}{2}}{\epsilon t-1}
\le m\cdot \forbmax(m-2,{\cal T}_{\ell}(3)\backslash {\cal T}_{\ell}(2)\cup  F).\label{A01bd}\end{E}
Using $\ncols{A}=2\cdot c(2,\ell)+\ncols{A_2}+\ncols{A_{01}}$, we obtain our desired bound.
\qed

\vskip 10pt

\section{Two-columned matrices}\label{twocolumn}

The main result of this section is the following. The proof is given after \lrf{underTk} and \lrf{underTk2}.
\begin{thm}\label{2colgeneral}
$\fb(m,3,\mathcal{T}_{\ell}(3)\backslash\mathcal{T}_{\ell}(2) \cup F)$ is $O\left(\sum^m_{k=0} \fbmax_k(m,F)\right)$ for all two-columned matrices $F$.
\end{thm}

We have some useful results for two-columned $F$.
Theorem 1.2 in \cite{AK}, gives us insight into $F_{0,b,b,0}$ with a strong stability result.
For our purposes we only need the following.

\begin{lemma} \cite{AK} Let $k,b$ be given with $b\ge 1$. Then $\forb_k(m,F_{0,b,b,0})$ is $O(m^{b-1})$. \label{k0bb0}\end{lemma}

Lemma 5.4 in \cite{AK},  repeated below,  gives us insight into $F_{0,b,b,1}$ which helps us consider $F_{1,b,b,1}$ for $b\ge 2$..

\begin{lemma}\label{k,1,b,b,1} Suppose $r \geq 1$ and 
$\mathcal{F}$ is a $k$-uniform family of subsets of $[m]$, with
$k \geq r+2$, so that every
pair $A,B\in \mathcal{F}$ is either disjoint or intersects
in at least $k-r$ points, and for every $A\in \mathcal{F}$ we have $1\notin A$. 
Then $|\mathcal{F}|$ is $O(m^r)$.
\end{lemma}

Translated in our language it says $\forb_k(m,F_{0,r+1,r+1,1})$ is $\Theta(m^r)$ for $k\ne r+1$.
Note that  $\forb_{r+1}(m,F_{0,r+1,r+1,1})$ is $\Theta(m^{r+1})$ by taking all columns of $r+1$ 1's.
By taking (0,1)-complements where $F_{0,r+1,r+1,1}^c=F_{1,r+1,r+1,0}$.  Thus for $k\ne m-r-1$, $\forb_k(m,F_{1,r+1,r+1,0})$ is $\Theta(m^r)$ and   $\forb_{r+1}(m,F_{1,r+1,r+1,0})$ is $\Theta(m^{r+1})$ by taking all columns of $r+1$ 1's.

\begin{cor} Let $r\ge 1$.  Let $m$ be given.  For $k\ne r+1,r+2,m-r-2,m-r-1$, we have $\forb_k(m,F_{1,r+1,r+1,1})$ is $\Theta(m^{r})$.
For $k= r+1,r+2,m-r-2$ or $m-r-1$ we have $\forb_k(m,F_{1,r+1,r+1,1})$ is $\Theta(m^{r+1})$.\label{k1bb1}\end{cor}
\proof Assume  $k\ne r+1,r+2,m-r-2,m-r-1$ and $\forb_k(m,F_{1,r+1,r+1,0})\le cm^{r}$. Let $A\in\Av(m,F_{1,r+1,r+1,1})$.  Consider  row 1. The number of 0's plus the number of 1's in row 1 is $\ncols{A}$.
 Let $B$ be the submatrix of $A$ formed by the columns with a 1 in row 1 and rows 
$2,3,\ldots m$. Then $B\in\Av_{k-1}(m-1,F_{0,r+1,r+1,1})$ and so $\ncols{B}\le\forb_{k-1}(m-1,F_{0,r+1,r+1,1})$. Note that $k-1\ne r+1$.
Let $C$ be the submatrix of $A$ formed by the columns with a 0 in row 1 and rows 
$2,3,\ldots m$. Then $C\in\Av_{k}(m-1,F_{1,r+1,r+1,0})$.  Now $F_{1,r+1,r+1,0}$ is the (0,1)-complement of $F_{0,r+1,r+1,1}$
and so $\forb_{k}(m-1,F_{1,r+1,r+1,0})=\forb_{m-k-1}(m-1,F_{0,r+1,r+1,1})$. Note that $k\ne (m-1)-r$. We deduce that
$\ncols{A}=\ncols{B}+\ncols{C}\le 2cm^{r}$ for $k\ne r+1,r+2,m-r-2,m-r-1$.  For the remaining cases  $k= r+1,r+2,m-r-2$ or $m-r-1$, we deduce that $\ncols{A}$ is $O(m^{r+1})$.
\qed 

\vskip 10pt
Let  a two-columned $F$ and an  $A \in \Av(m,3,\mathcal{T}_{\ell}(3)\backslash\mathcal{T}_{\ell}(2) \cup F)$ with fixed column sum $k$ be given. Note that by \trf{columnsumkextra}, given $L$, there exists a constant $t$ such that after $t\left(\sum_{i=k-t}^k \forbmax_i(m,F)\right)$ columns, we either find in $A$ a configuration of $\mathcal{T}_L(3)\backslash\mathcal{T}_L(2)$, or $T_L(0,2,1)$ or $F$.
 If $L \geq \ell$, the only object on this list not forbidden is $T_L(0,2,1)$, so we may assume we find this configuration. Note that $T_{L/2}(1,0) \prec T_L(0,2,1)$, so $T_{L/2}(1,0)$ must appear. Reorder the columns so that the 1's are above the diagonal in $T_L(0,1)$. Now, using the previous notation for two-columned matrices, let $F = F_{a,b,c,d}$. Notice that if we delete the first $a$ columns of $T_L(0,1)$, every pair of columns has $a$ copies of $[1\,1]$; if we delete the last $d$ columns, every pair of columns has $d$ copies of $[0\,0]$; and if we take every $c$th column of what remains, every pair of columns has  $c$ copies of $[0\,1]$. \\

Let $A'$ be the submatrix of $A$ obtained by taking the selected columns from $T_L(0,1)$ and deleting the rows from $T_L(0,1)$. Note that in the deleted rows, every pair of columns has $a$ copies of $[1\,1]$, $c$ copies of $[0\,1]$, and $d$ copies of $[0\,0]$, so if any pair of columns of $A'$ have $b$ copies of $[1\,0]$, $A$ contains $F$. Also, since $A$ has fixed column sum, and the column sums of $T_L(0,1)$ increase from left to right, the column sums of $A'$ decrease from left to right. To use these facts we need the following lemma.

\begin{lemma}\label{underTk} Let $M$ be an $m$-rowed matrix such that:
\begin{enumerate}
\item[(i)] $M$ does not contain $\begin{bmatrix} \mathbf{1}_b & \mathbf{0}_b \end{bmatrix}$ as a $b\times 2$ submatrix for some $b$.
\item[(ii)] If $i < j$, column $i$ of $M$ has more 1's than column $j$
\item[(iii)] $M$ avoids $\mathcal{T}_{\ell}(3)\backslash\mathcal{T}_{\ell}(2)$
\end{enumerate}
Then for every $r$, there exists a constant $c_r$ (dependent on $k$) such that if $M$ has more than $c_r$ columns, it has an $r \times r$ configuration in ${\cal P}_r^2$ \rf{Ptemplate} with 1's on the diagonal and 2's below the diagonal.
\end{lemma}

\proof
We proceed by induction. When $r = 1$, the desired object is just a single 0, so the lemma is trivial. Suppose the lemma holds for $r$. We claim that the lemma holds for $r + 1$ with $c_{r+1} = R_5(\ell,\ell,\ell,c_r,b + 1) + b$. Suppose $M$ satisfies the hypotheses of the lemma and has $c_{r+1}$ columns. Define $M'$ to be the restriction of $M$ to the rows with a 1 in the first column. Since the column sums of $M$ strictly decrease from left to right, the $(b + 1)$th column of $M$ has at least $b$ fewer 1's than the first, which implies that there must be at least $b$ non-1 entries in the $(b + 1)$th column of $M'$. At most $b - 1$ of these entries are 0 by condition (i), so there is at least one 2. Pick one. The $(b + 2)$th column of $M'$ has at least two 2's, at least one of which is in a different row than the one already chosen. Pick one such 2. Similarly the $(b + 3)$th column of $M'$ has a 2 in a different row than the 2's already selected, and so on; continuing in this way, we find a diagonal of 2's of length $\ncols{M} - b = R_5(\ell,\ell,\ell,c_r,b + 1)$. Let the square submatrix of $A$ induced by the row and column indices of the chosen diagonal be $M''$. \\

We now produce a colouring of the complete graph on $\ncols{M''}$ vertices as follows. Given $i < j$, if $M''_{ij},M''_{ji} \neq 0$, colour edge $\{i,j\}$ with the ordered pair $(M''_{ij},M''_{ji})$; if $M''_{ij} = 0$ or $M''_{ji} = 0$, colour $\{i,j\}$ with 0. Now there are five colours: $(1,1),(1,2),(2,1),(2,2)$, and $0$. By Ramsey Theory, we have a clique of size $\ell$ of colour $(1,1),(1,2)$ or $(2,1)$ or a clique of size $b + 1$ of colour 0 or a clique of size $c_r$ of colour $(2,2)$. In the first case, all three colours give rise to a member of $\mathcal{T}_{\ell}(3)\backslash\mathcal{T}_{\ell}(2)$. In the second case, we have a column with $b$ 0's opposite the 1's in the first column of $M$, contradicting condition (i). Hence the only allowed case is the third, which corresponds to a block of 2's. In particular, there is a row $x$ of $M$ with $c_r$ 2's and a 1. Look under the 2's; the resulting matrix has $c_r$ columns and satisfy the hypotheses of the lemma, so by induction there is an $r \times r$ configuration with 1's on the diagonal and 2's above. Adding in row $x$ gives an $(r + 1) \times (r + 1)$ configuration of the desired type.
\qed

\begin{lemma}\label{underTk2}
Let $M$ satisfy the hypotheses of Lemma \ref{underTk} with $c_r$ defined there. Then $\ncols{M} < c_{R_3(2\ell,\ell,\ell)}$.
\end{lemma}

\proof
We use the notation $c_r$  from the statement of \lrf{underTk}.  Suppose $\ncols{M}\ge c_{R_3(2\ell,\ell,\ell)}$. By \lrf{underTk}, there exists a configuration $N\prec M$ in ${\cal P}^2_{R_3(2\ell,\ell,\ell)}$ 
 with 1's on the diagonal and 2's below the diagonal. We colour a complete graph $K_{R_3(2\ell,\ell,\ell)}$ as follows: for $i < j$, colour edge $(i,j)$ with $N_{ji}$ (note that $N_{ij} = 2$). By the definition of $R_3(2\ell,\ell,\ell)$, there is a monochromatic clique. Three cases are possible. If there is a clique of size $2\ell$ of colour 0, we get $T_{2\ell}(0,1,2)$, which contains $T_{\ell}(0,2)$. If there is a clique of size $\ell$ of colour 1, we have $T_{\ell}(1,2)$, and a clique of size $\ell$ of colour 2 yields $I_{\ell}(1,2)$. This contradicts our assumption that $\ncols{M} \geq c_{R(2\ell,\ell,\ell)}$.
\qed
\vskip 10pt
\noindent{\bf Proof of \trf{2colgeneral}}:
Let $A \in \Av(m,3,\mathcal{T}_{\ell}(3)\backslash\mathcal{T}_{\ell}(2) \cup F)$ be given, with fixed column sum $k$. By \trf{columnsumkextra}, there exist constants $C$ and $d$ independent of $k$ such with more than $C(\sum^k_{i=k-d} \forbmax_i(m,F))$ columns, we either have one of the forbidden objects or a very large triangular matrix $T_t(0,2,1)$ with $t=c_{R_3(2\ell,\ell,\ell)}$. This yields a matrix $M$ satisfying the hypotheses of \lrf{underTk} with $\ncols{M}>c_{R_3(2\ell,\ell,\ell)}$. Then \lrf{underTk2} yields a contradiction.   Hence, $||A|| \leq C(\sum^k_{i=k-d} \forbmax_i(m,F))$, so $\forbmax_k(m,3,\mathcal{T}_{\ell}(3)\backslash\mathcal{T}_{\ell}(2) \cup F) \leq C(\sum^k_{i=k-d} \forbmax_i(m,F))$. Summing over $k$ gives the desired result.
\qed
\vskip 10pt
This result can be used to give bounds for many 2-columned matrices. \\

\noindent{\bf Proof of Theorem \ref{abba}}:
For $F_{0,b,b,0}$ we use \lrf{k0bb0} which yields that \hfil\break$\forbmax_k(m,3,\mathcal{T}_{\ell}(3)\backslash\mathcal{T}_{\ell}(2) \cup F_{0,b,b,0})$ is $O(m^{b-1})$. Then \trf{2colgeneral} yields that \hfil\break$\forb(m,3,\mathcal{T}_{\ell}(3)\backslash\mathcal{T}_{\ell}(2) \cup F_{0,b,b,0})$ is $O(m^b)$.  From \cite{AK}, $\forb(m,F_{0,b,b,0})$ is $\Theta(m^b)$.

By \crf{k1bb1}, $\forbmax_k(m,F_{1,b,b,1})$ is $O(m^{b-1})$ for $b \geq 2$ and $k\ne r+1,r+2,m-r-2,m-r-1$.  For $k=r+1,r+2,m-r-2$ or $m-r-1$, the bound os $O(m^b)$.  By \trf{2colgeneral}, $\forb(m,3,\mathcal{T}_{\ell}(3)\backslash\mathcal{T}_{\ell}(2) \cup F_{1,b,b,1})$ is $O(m^b)$. From \trf{rows0and1} we may extend this to obtain $\forb(m,3,\mathcal{T}_{\ell}(3)\backslash\mathcal{T}_{\ell}(2) \cup F_{a,b,b,a})$ is $O(m^{a + b - 1})$. This is the correct bound by \cite{AK}. 
\qed
\vskip 10pt
\noindent{\bf Proof of Theorem \ref{I2}}: Use \lrf{k0bb0} with $b=1$ which by \trf{2colgeneral}, yields that $\forb(m,3,\mathcal{T}_{\ell}(3)\backslash\mathcal{T}_{\ell}(2) \cup F_{0,1,1,0})$ is $O(m)$. 
\qed
\vskip 10pt
We do not know how to do solve for $F=F_{1,1,1,1}$ for which $\forb(m,F_{1,1,1,1})$ and $\forb_k(m,F_{1,1,1,1})$ are both $\Theta(m)$. 
Similarly, the case $F=F_{a,1,1,a}$ for $a\ge 2$ is not solved. We have $\forb(m,F_{a,1,1,a})$ is $\Omega(m^{a})$.
The following results give bounds which must be close to the correct bounds.

\begin{thm}\label{1,1,1,1}
$\forb(m,3,\mathcal{T}_{\ell}(3)\backslash\mathcal{T}_{\ell}(2) \cup F_{1,1,1,1})$ is $O(m\log m)$
\end{thm}
\proof
Let $A \in \Av(m,F_{1,1,1,1})$ with column sum $k$ be given. If $k \leq m/2$, then every pair of columns has a $[0\,0]$. Since the column sum is fixed, every pair of columns has an $I_2$. Hence there must be no $[1\,1]$ in any pair of columns. This means the 1's must all appear on disjoint rows, so there are at most $\frac{m}{k}$ columns. If $k > m/2$, take the 0-1 complement to get a similar result. Applying \trf{2colgeneral} and summing over $k$ gives the result.  
\qed
\vskip 10pt
Applying \trf{rows0and1} gives the following corollary.

\begin{cor}
$\forb(m,3,\mathcal{T}_{\ell}(3)\backslash\mathcal{T}_{\ell}(2) \cup F_{a,1,1,a})$ is $O(m^a\log m)$
\end{cor}

Of course, the extra factor of $\log m$ is undesirable. However, given that all known forbidden families have a polynomial bound, this strongly suggests that the actual bound for $F_{a,1,1,a}$ is $O(m^a)$. \\

\section{An example with 3 columns}

Define the useful notation $A|_S$ to denote the submatrix of $A$ given by the rows $S$.
In order to prove \trf{linear} for $H$ given in \rf{H}, we find the following lemma usful. A standard decomposition applied to 3-matrices considers 
 deleting a row $i$ from  a simple 3-matrix $A$. The resulting matrix might not be simple. Let $C_{a,b}(i)$ be the simple 3-matrix that consists of the repeated columns of the matrix that is obtained when deleting row $r$ from $A$ that lie under both symbol $a$ and $b$ in row $i$. In particular $[a\,b]\times C_{a,b}(i)\prec A$.   Let $B(i)$ denote the $(m-1)$-rowed simple 3-matrix obtained from $A$ by deleting row $i$ and any repeats of  columns so that 
 $$\ncols{A}\le\ncols{B(i)}+\ncols{C_{0,1}(i)}+\ncols{C_{1,2}(i)}+\ncols{C_{0,2}(i)}.$$
  The inequality arises from columns that are repeated three times in the matrix obtained from $A$ by deleting row $r$ but get counted four times on the right hand side. This bound on $\ncols{A}$  is often amenable to induction on the number of rows. 
If $K_2=[0\,1]\times [0\,1]\not\prec A$,  or in our case $H\not\prec A$, we  deduce that $\ncols{C_{0,1}(i)}$ is $O(1)$, namely the constant bound for 
$\forb(m,3,{\cal T}_{\ell}(3)\backslash {\cal T}_{\ell}(2)\cup [0\,1])$ using \crf{[01]}.  The following lemma could help with 
$\ncols{C_{1,2}(i)}$ and $\ncols{C_{0,2}(i)}$.

\begin{lemma}Let $A\in\Av(m,3,{\cal T}_{\ell}(3)\backslash {\cal T}_{\ell}(2))$. Assume 
for some set of rows $S$  we have $[\0\,|\,I_{|S|}]\prec A|_S$ and for each pair of rows $i,j\in S$, we have no
$[\linelessfrac{1}{1}]$ in $A$.   If $|S|>3\ell\cdot c(2,\ell)$, then there is some row $i\in S$
for which $\ncols{C_{1,2}(i)}=\ncols{C_{0,2}(i)}=0$. \label{identity}\end{lemma}
\proof 
We will show that $\ncols{C_{1,2}(i)}>0$ for only a few choices $i\in S$ and similarly show that $\ncols{C_{0,2}(i)}>0$ for only a few choices $i\in S$.   Then for $S$ large enough, there will be some $i\in S$ with $\ncols{C_{1,2}(i)}=\ncols{C_{0,2}(i)}=0$.

  Let $U$ denote the rows $i\in S$ for which  $\ncols{C_{1,2}(i)}>0$. Assume $|U|\ge \ell\cdot c(2,\ell)$.  
  When $\ncols{C_{1,2}(i)}>0$, we have (at least) two columns in $A$ differing only in row $i$, one with a 1 and one with a 2. Choose one such pair of columns $\gamma$,$\delta$ as shown:
$$\begin{array}{c}i\\ U\backslash i \\ {[m]\backslash U}\\  \end{array}
\left[\begin{array}{cc}1&2\\ \alpha&\alpha\\ \beta&\beta\\ \end{array}\right]\prec A.$$
It is possible that for many $i$, the same second column might be chosen.
By the property of $A$ that $A$ has no $\left[\linelessfrac{1}{1}\right]$ on rows of $S$ and hence $U$, we deduce   
$\delta|_U=[\linelessfrac{2}{\alpha}]$ is a (0,2)-vector.  By \trf{constantbd} (in this case due to \cite{BaBo}), we have that there are at most $c(2,\ell)$ choices for 
$[\linelessfrac{2}{\alpha}]$. Now there are $|U|$ choices for $i$ and so, given our bound on $|U|$,  there are $\ell$ choices for $i\in U$ which have the same $[\linelessfrac{2}{\alpha}]$.  Now considering the $\ell$ columns $[\linelessfrac{1}{\alpha}]$ yields an  $\ell\times\ell$ matrix in $A|_U$ with 1's on diagonal and 2's off the diagonal namely $I_{\ell}(2,1)\in {\cal T}_{\ell}(3)\backslash {\cal T}_{\ell}(2)$, a contradiction.  Thus $\ncols{C_{1,2}(i)}>0$ for less than $\ell\cdot c(2,\ell)$ choices $i$. 

Assume $\ncols{C_{0,2}(i)}>0$ for $2\ell\cdot c(2,\ell)$ choices $i$.  Denote the choices by $V$. Then we have the following
\begin{E}\begin{array}{c}i\\ V\backslash i\\  {[m]\backslash V}\\  \end{array}\left[\begin{array}{cc}0&2\\ \alpha&\alpha\\ \beta&\beta\\ \end{array}\right]\prec A\label{pair}\end{E}
This case is a little more complicated because $\alpha$ may have up to one 1.  
 We choose a subset $W\subseteq V$ of the rows where $\alpha$ has no 1's.
This can be done as follows. Choose some row $i_1\in V$ and assume the corresponding  choice of columns yields an $\alpha$ with  a 1 in row $j_1\in V$  and if not let $j_1=i_1$.  Now choose a row $i_2\in V\backslash\{i_1,j_1\}$ and assume the corresponding  $\alpha$ has a 1 in row $j_2\in V$ and if not $j_2=i_2$. Now choose a row $i_3\in V\backslash\{i_1,j_1,i_2,j_2\}$ and assume the corresponding  $\alpha$ has a 1 in row $j_3\in V$ and if not $j_3=i_3$.
Continue in this way to form $W=\{i_1,i_2,\ldots ,i_{\ell\cdot c(2,\ell)}\}$ using the fact $|V| \ge 2\ell\cdot c(2,\ell)$. 
$|W|\ge \ell\cdot c(2,\ell))$ and for each $i\in W$ we have $\ncols{C_{0,2}(i)}>0$ where we have on pair of cols in $A$ as in \rf{pair} with $\alpha$ having no 1's.  Now repeat the above argument for the (1,2)-case to obtain an $\ell\times \ell$ matrix in $A$    with 0's on diagonal and 2's off the diagonal, namely $I_{\ell}(2,0) \in {\cal T}_{\ell}(3)\backslash {\cal T}_{\ell}(2)$, a contradiction.   
 Thus $\ncols{C_{0,2}(i)}>0$ for less than $2\ell\cdot c(2,\ell)$ choices $i$.

We deduce that for $|S|>3\ell\cdot c(2,\ell)$, there exists a row $i$ with $|C_{1,2}(i)|=|C_{0,2}(i)|=0$.
\qed 

\vskip 10pt

\noindent {\bf Proof of Theorem \ref{linear}}:
Let $A \in \Av(m,3,\mathcal{T}_{\ell}(3)\backslash\mathcal{T}_{\ell}(2) \cup H)$ be given. If $A$ contains a large identity (or its complement) then by \lrf{identity} there exists a row $i$ with $C_{1,2}(i) = C_{0,2}(i) = \emptyset$. Note that $\ncols{C_{0,1}(i)} $ is $O(1)$ by \crf{[01]}, since $C_{0,1}(i)$ avoids $[0\,1]$. Thus, we can delete row $i$ and at most $O(1)$ columns and obtain a simple matrix.  Then induction on $m$ would yield the desired $O(m)$ bound.  Our goal is to show that a large identity must occur. \\

Let $A_k$ be the submatrix of $A$ with column sum $k$. If, for any $L$, $\ncols{A_k}>c(3,L)$ then either $I_L(0,1)\prec A_k$, $I_L(1,0)\prec A_k$, or $T_L(0,1)\prec A_k$. We will take $L$ to be large.   We note that $H\prec I_L(0,1)$ and so the first case does not occur. In the second case $I_L(1,0)\prec A_k$ we have that $A_k$ and indeed $A$ does not have $[\linelessfrac{0}{0}]$ on the $L$ rows containing $I_L(1,0)$ else $H\prec A$.  Then  apply \lrf{identity}, using the (0,1)-complement and  note that 
$\ncols{C_{1,0}(i)}$ is $O(1)$ by \crf{[01]} since $H\prec [0\,1]\times [0\,1]$ and hence $[0\,1]\not\prec C_{0,1}(i)$. This yields that $\ncols{A}$ is $O(m)$ by induction on $m$. 

   In the third case, with the triangular matrix $T_L(0,1)\prec A_k$, let $A_k'$ be the matrix consisting of the columns from $A_k$ 
    containing $T_L(0,1)$. Assume the columns of $A_k'$ are ordered consistent with $T_L(0,1)$. Let $A_k''$ be the submatrix obtained from $A_k'$ by deleting the $L$ rows containing $T_L(0,1)$.  Then  the column sums of $A_k''$ are decreasing from left to right.  Let $S$ be the rows containing 1's in the first column of $A_k''$.    Every triple of columns in $T_L(0,1)$ has the submatrix $[0\,0\,1]$, so $A_k''$ does not contain any submatrix $[1\,0\,0]$ else $H\prec A_k''\prec A$. Thus $A_k''|_S$ does not contain $[0\,0]$. Also $A_k''$ has decreasing column sums from left to right. We proceed in a manner similar to the proof of \lrf{underTk}.  
    We first find a diagonal of  entries either 0 or 2. By the pigeonhole principle, there is a long diagonal of 2's or a long diagonal of 0's. If there is a long diagonal of 2's we apply Ramsey Theory as before. Large cliques involving 0's are not allowed since $[0\,0]$ is forbidden, and hence we are forced to have a block of 2's. This yields a submatrix $[1\,2\,2\cdots 2]$ and so we proceed, as in proof of \lrf{underTk}, considering the columns containing the 2's. If we continue doing this for enough rows, we find a forbidden object. Hence, there must be a point where no sufficiently long diagonal of 2's exists, so there is a long diagonal of 0's. In this case, apply Ramsey Theory again. Recalling that we have no submatrix $[0\,0]$, the only configurations that result are either in ${\cal T}_{\ell}(3)\backslash {\cal T}_{\ell}(2)$ or an identity complement $I_t(1,0)$  for some large $t$.  
    Given that $H\not\prec A$, we have that there is no $\left[\linelessfrac{0}{0}\right]$ on any pair of the $t$ rows. This  allows us to use \lrf{identity}. 

If $\ncols{A_k}$ is bounded by a constant for all $k$ then $\ncols{A}$ is $O(m)$. If $\ncols{A_k}$ is a big enough constant  then we obtain an $I_t(1,0)\prec A_k$  for some appropriately large $t$. By \lrf{identity} we find some $i\in [m]$ with  $\ncols{C_{1,2}(i)}=\ncols{C_{0,2}(i)}=0$.   As noted above, 
$\ncols{C_{0,1}(i)}$ is $O(1)$.   Thus we can delete  row $i$ of $A$ and at most $O(1)$ columns from $A$ to obtain a simple matrix in 
$\Av(m-1,3,{\cal T}_{\ell}(3)\backslash {\cal T}_{\ell}(2)\cup H)$ and then apply induction. 
\qed

 \section{Open problems}
 
 Some small examples of $F$ for which we have not handled 
 $\forb(m,3,{\cal T}_{\ell}(3)\backslash {\cal T}_{\ell}(2)\cup F)$ include:
 
 $$K_2=\left[\begin{array}{cccc}0&0&1&1\\ 0&1&0&1\\ \end{array}\right], \hbox{ bound should be }O(m)$$
 $$F_{0,2,1,0}=\left[\begin{array}{cccc}1&0\\1&0\\ 0&1\\  \end{array}\right], \hbox{ bound should be }O(m)$$
 $$F_{1,1,1,0}=\left[\begin{array}{cccc}1&1\\1&0\\ 0&1\\  \end{array}\right], \hbox{ bound should be }O(m)$$
 We would particularly like to have a general result that $\forb(m,3,{\cal T}_{\ell}(3)\backslash {\cal T}_{\ell}(2)\cup ([0\,1]\times F))$ is 
 $O(m\times\forb(m,3,{\cal T}_{\ell}(3)\backslash {\cal T}_{\ell}(2)\cup  F))$ matching the standard induction results for (0,1)-forbidden configurations.
 
 Given a $(0,1)$-column $\alpha$,  we might consider a 3-matrix $A\in\Av(m, 3,{\cal T}_{\ell}(3)\backslash {\cal T}_{\ell}(2))$ such that each column of $A$ arises from $\alpha$ by setting certain entries to 2.  We deduce that $[0\,1]\not\prec A$ and so by \crf{[01]}, we have the  interesting fact that $\ncols{A}$ is $O(1)$.  In some sense the columns of $A$ are a 3-matrix replacement for $\alpha$.  We were unable to exploit this for \probrf{problem}.

 \end{document}